\documentclass[11pt]{amsart}

\usepackage[a4paper,hmargin=3.5cm,vmargin=4cm]{geometry}
\usepackage{amsfonts,amssymb,amscd,amstext}
\usepackage{graphicx}
\usepackage[dvips]{epsfig}
\usepackage{pstricks,pst-plot}
\usepackage{float}

\usepackage{fancyhdr}
\input xy
\xyoption{all}

\usepackage{times}

\usepackage{enumerate}
\usepackage{titlesec}
\usepackage{mathrsfs}

\titleformat{\section}
{\filcenter\bfseries\large} {\thesection{.}}{0.2cm}{}
\titleformat{\subsection}[runin]
{\bfseries} {\thesubsection{.}}{0.15cm}{}[.]
\titleformat{\subsubsection}[runin]
{\em}{\thesubsubsection{.}}{0.15cm}{}[.]

\usepackage[up,bf]{caption}

\newtheorem{theorem}{Theorem}[section]

\newtheorem{lemma}[theorem]{Lemma}
\newtheorem{corollary}[theorem]{Corollary}

\theoremstyle{definition}

\numberwithin{equation}{section}
\numberwithin{figure}{section}

\newcommand\Rcal{\mathcal{R}}

\newcommand\Cscr{\mathscr{C}}

\newcommand\Lscr{\mathscr{L}}
\newcommand\Oscr{\mathscr{O}}

\def\c{\mathbb{C}}

\renewcommand\d{\mathbb D}

\def\n{\mathbb{N}}
\renewcommand\r{\mathbb{R}}

\newcommand\z{\mathbb{Z}}

\newcommand\igot{\mathfrak{i}}

\renewcommand\igot{\mathfrak{i}}

\newcommand\pgot{\mathfrak{p}}

\newcommand\Agot{\mathfrak{A}}
\newcommand\Igot{\mathfrak{I}}

\renewcommand\imath{\igot}

\newcommand\cM{\overline{M}}

\newcommand\wt{\widetilde}

\newcommand\di{\partial}

\newcommand\dist{\mathrm{dist}}

\newcommand\length{\mathrm{length}}

\newcommand\Flux{\mathrm{Flux}}

\newcommand\CMI{\mathrm{CMI}}

\usepackage{color}

\begin{document}

\fancyhead[LO]{Complete dense minimal surfaces}
\fancyhead[RE]{A.\ Alarc\'on and I.\ Castro-Infantes}
\fancyhead[RO,LE]{\thepage}

\thispagestyle{empty}

\vspace*{7mm}
\begin{center}
{\bf \LARGE Complete minimal surfaces densely lying in

arbitrary domains of $\r^n$}
\vspace*{5mm}

{\large\bf Antonio Alarc\'on \; and \; Ildefonso Castro-Infantes}
\end{center}

\vspace*{7mm}

\begin{quote}
{\small
\noindent {\bf Abstract}\hspace*{0.1cm}
In this paper we prove that, given an open Riemann surface $M$ and an integer $n\ge 3$, the set of complete conformal minimal immersions $M\to\r^n$ with $\overline{X(M)}=\r^n$ forms a dense subset in the space of all conformal minimal immersions $M\to\r^n$ endowed with the compact-open topology.
Moreover, we show that every domain in $\r^n$ contains complete minimal surfaces which are dense on it and have arbitrary orientable topology (possibly infinite); we also provide such surfaces whose complex structure is any given bordered Riemann surface.

Our method of proof can be adapted to give analogous results for non-orientable minimal surfaces in $\r^n$ $(n\ge 3)$, complex curves in $\c^n$ $(n\ge 2)$, holomorphic null curves in $\c^n$ $(n\ge 3)$, and holomorphic Legendrian curves in $\c^{2n+1}$ $(n\in\n)$.

\vspace*{0.1cm}
\noindent{\bf Keywords}\hspace*{0.1cm} Complete minimal surfaces, Riemann surfaces, holomorphic curves.

\vspace*{0.1cm}

\noindent{\bf MSC (2010):}\hspace*{0.1cm} 49Q05, 32H02}
\end{quote}


\section{Introduction and main results}\label{sec:intro}

The existence of complete minimal surfaces densely lying in $\r^3$ is well-known. 
The first example of such, due to Rosenberg, was obtained by Schwarzian reflection on a fundamental domain, 
is simply-conneted, and has bounded curvature. Later, G\'alvez and Mira \cite{GalvezMira2004BBMS} found complete dense simply-connected minimal surfaces in $\r^3$, in explicit coordinates, as solution to certain Bj\"orling problems. Finally, L\'opez \cite{Lopez2014JGA} constructed complete dense minimal surfaces in $\r^3$ with weak finite total curvature, arbitrary genus, 
and parabolic conformal type; so far, these are the only known examples with non-trivial topology. 
In a parallel line of results, Andrade \cite{Andrade2000PAMS} gave an example of a complete simply-connected minimal surface in $\r^3$ which is not dense in the whole space but its closure has nonempty interior. It is therefore a natural question whether a given domain in $\r^3$ 
contains complete minimal surfaces which are dense on it; as far as the authors knowledge extends, no domain is known to enjoy this property besides $\r^3$ itself.

The aim of this paper is to answer the above question by showing a general existence result for complete dense minimal surfaces in {\em any} given domain $D\subset\r^n$ for arbitrary dimension $n\ge 3$. We provide such surfaces with {\em arbitrary orientable topology and flux map}; moreover, if $n\ge 5$ we give examples with no self-intersections. Furthermore, if $D=\r^n$ then we construct such surfaces not only with arbitrary topology but also with {\em arbitrary complex structure}. To be precise, our first main result may be stated as follows.

\begin{theorem}\label{th:intro-main-v2}
Let $D\subset\r^n$ $(n\ge 3)$ be a domain, $M$ be an open Riemann surface, $\pgot\colon H_1(M;\z)\to\r^n$ be a group homomorphism, $K\subset M$ be a smoothly bounded Runge compact domain, and $X\colon K\to\r^n$ be a conformal minimal immersion of class $\Cscr^1(K)$. Assume that $X(K)\subset D$ and that the flux map $\Flux_X\colon H_1(K;\z)\to\r^n$ of $X$ satisfies $\Flux_X(\gamma)=\pgot(\gamma)$ for all closed curves $\gamma\subset K$.

Then, for any $\epsilon>0$, there are a domain $\Omega\subset M$ and a complete conformal minimal immersion $Y\colon \Omega\to \r^n$ satisfying the following properties:
\begin{enumerate}[\rm (I)]
\item $K\subset \Omega$ and $\Omega$ is a deformation retract of $M$ and homeomorphic to $M$.
\item $\|Y-X\|_{1,K}<\epsilon$.
\item $\Flux_Y(\gamma)=\pgot(\gamma)$ for all closed curves $\gamma\subset\Omega$.
\item $Y(\Omega)\subset D$ and the closure $\overline{Y(\Omega)}= \overline{D}$.
\item $Y$ is one-to-one if $n\ge 5$.
\end{enumerate}
Furthermore, if $D=\r^n$ we may choose $\Omega=M$. 
\end{theorem}

Theorem \ref{th:intro-main-v2} gives the first examples of complete dense minimal surfaces in $\r^n$ for $n>3$. Notice that the density of $Y(M)$ in $D$ does not allow the immersions $Y\colon\Omega\to D$ in the theorem to be proper maps.

We emphasize that, although certainly wild, complete dense minimal surfaces in $\r^n$ $(n\ge 3)$ are surprisingly abundant. Indeed, if we denote by $\CMI(M,\r^n)$ the space of all conformal minimal immersions of a given open Riemann surface $M$ into $\r^n$ (which is nonempty by the results in Alarc\'on and L\'opez \cite{AlarconLopez2012JDG}), Theorem \ref{th:intro-main-v2} ensures that 
{\em those conformal minimal immersions $M\to\r^n$ which are complete and have dense image form a dense subset of $\CMI(M,\r^n)$ with respect to the compact-open topology}.

It is also worth mentioning at this point that it is not hard to find dense minimal surfaces in $\r^n$ for any $n\ge 3$. Indeed, solving the Bj\"orling problem for any real analytic regular dense curve in $\mathbb{R}^n$ 
and any tangent plane distribution along it gives such a surface; we thank Pablo Mira for providing us with this simple argument. Obviously, this method only produces simply-connected examples and does not guarantee their completeness. As will become apparent later in this introduction, constructing {\em complete} dense minimal surfaces in $\r^n$, {\em prescribing their topology and even their complex structure}, is a much more arduous task which requires of a number of powerful and sophisticated tools of the theory that have been developed only recently.

It is well-known that a general domain $D\subset\r^n$ does not contain minimal surfaces with arbitrary complex structure. Indeed, if for instance $D$ is relatively compact then it only admits minimal surfaces of {\em hyperbolic} conformal type (see \cite{FarkasKrabook}). We also prove in this paper that every domain $D\subset\r^n$ contains complete minimal surfaces which are dense on it and whose complex structure is any given bordered Riemann surface. 
\begin{theorem}\label{th:intro-main}
Let $D\subset\r^n$ $(n\ge 3)$ be a domain and $\cM=M\cup bM$ be a compact bordered Riemann surface. Every conformal minimal immersion $X\colon \cM\to \r^n$ of class $\Cscr^1(\cM)$, with $X(\cM)\subset D$, may be approximated uniformly on compact subsets of $M=\cM\setminus bM$ by complete conformal minimal immersions $Y\colon M\to \r^n$ assuming values in $D$ and such that $\overline{Y(M)}=\overline D$ and $\Flux_Y=\Flux_X$. Moreover, if $n\ge 5$ then the approximating immersions $Y$ can be chosen to be one-to-one.
\end{theorem}

Recall that a {\em compact bordered Riemann surface} is a compact Riemann surface $\cM$ with nonempty boundary $bM\subset\overline M$ consisting of finitely many pairwise disjoint smooth Jordan curves. The interior $M=\cM\setminus bM$ of $\cM$ is called a {\em bordered Riemann surface}. 
By a {\em conformal minimal immersion $\overline M\to\r^n$ of class $\Cscr^1(\overline M)$} we mean a map of class $\Cscr^1(\overline M)$ whose restriction to $M$ is a conformal minimal immersion.

%
%
We shall prove Theorems \ref{th:intro-main-v2} and \ref{th:intro-main} in Section \ref{sec:MR}. 
The main tools in our method of proof come from the strong connection between minimal surfaces in $\r^n$ and Complex Analysis; in particular, {\em Oka theory} (see the note by L\'arusson \cite{Larusson2010NAMS} and the surveys by Forstneri\v c and L\'arusson \cite{ForstnericLarusson2011NY}, Forstneri\v c \cite{Forstneric2013AFSTM}, and Kutzschebauch \cite{Kutzschebauch2014SPMS} for an introduction to this theory, and the monograph by Forstneri\v c \cite{Forstneric2011-book} for a comprehensive treatment; see e.g. Alarc\'on and Forstneri\v c \cite{AlarconForstneric2014IM,AlarconForstneric2015AS} or Alarc\'on, Forstneri\v c, and L\'opez \cite{AlarconForstnericLopez2016MZ}, and the references therein, for a discussion of the interplay between minimal surfaces and Oka manifolds). To be more precise, our proof relies on a {\em Runge-Mergelyan type approximation theorem} for conformal minimal immersions of open Riemann surfaces into $\r^n$ (see Alarc\'on and L\'opez \cite{AlarconLopez2012JDG} for $n=3$ and Alarc\'on, Forstneri\v c, and L\'opez \cite{AlarconForstnericLopez2016MZ} for arbitrary dimension), a {\em general position theorem} for conformal minimal surfaces in $\r^n$ for $n\ge 5$ (see \cite{AlarconForstnericLopez2016MZ}), and the existence of approximate solutions to certain {\em Riemann-Hilbert type boundary value problems} for conformal minimal surfaces in $\r^n$ where the complex structure of the central surface is a compact bordered Riemann surface (see Alarc\'on and Forstneri\v c \cite{AlarconForstneric2015MA} for $n=3$ and Alarc\'on, Drinovec Drnov\v sek, Forstneri\v c, and L\'opez \cite{AlarconDrinovecForstnericLopez2015PLMS} for $n\ge 3$). Actually, the Riemann-Hilbert method is not explicitly applied in the present paper but it plays a fundamental role in the proof of \cite[Lemma 4.1]{AlarconDrinovecForstnericLopez2015PLMS}, which we use in a strong way. Furthermore, our method of proof also exploits the technique by Forstneri\v c and Wold \cite{ForstnericWold2009JMPA} for {\em exposing boundary points on a bordered Riemann surface}, which pertains to Riemann Surface Theory. 

%
%
All the above mentioned tools are also available for some other families of surfaces which are the focus of interest, namely, {\em non-orientable minimal surfaces} in $\r^n$ for $n\ge 3$, {\em complex curves} in the complex Euclidean spaces $\c^n$ for $n\ge 2$, {\em holomorphic null curves} in $\c^n$ for $n\ge 3$, and {\em holomorphic Legendrian curves} in $\c^{2n+1}$ for $n\in\n$.
Thus, our methods easily adapt to give results analogous to Theorems \ref{th:intro-main-v2} and \ref{th:intro-main} in all these geometric contexts; we motivate, state, and discuss some of them in Section \ref{sec:results}.


\section{Preliminaries}\label{sec:prelim}

Given $n\in\n=\{1,2,3,\ldots\}$, we denote by $|\cdot|$, $\dist(\cdot,\cdot)$, and $\length(\cdot)$
the Euclidean norm, distance, and length
in $\r^n$, respectively. Given a set $A\subset \r^n$ we denote by $\overline A$ the topological closure of $A$ in $\r^n$. 

If $K$ is a compact topological space and $f\colon K\to \r^n$ is a continuous map, we denote by 
\[
\|f\|_{0,K}:=\max\{|f(p)|\colon p\in K\}
\]
the maximum norm of $f$ on $K$. 
If $K$ is a subset of a Riemann surface $M$, then for any $r\in\z_+=\n\cup\{0\}$ we denote by 
\[
	\|f\|_{r,K}
\] 
the standard $\Cscr^r$ norm of a function $f\colon K\to\r^n$ of class $\Cscr^r(K)$, where the derivatives are measured with respect to a fixed Riemannian metric on $M$ (the precise choice of the metric will not be important).

Given a smooth connected surface $S$ (possibly with nonempty boundary) and a smooth immersion $X\colon S\to\r^n$ $(n\ge 3)$, we denote by 
\[
	\dist_X\colon S\times S\to\r_+=[0,+\infty)
\]
the Riemannian distance induced on $S$ by the Euclidean metric of $\r^n$ via $X$:
\[
	\dist_X(p,q)=\inf\{\length(X(\gamma))\colon 
	\text{$\gamma\subset S$ arc connecting $p$ and $q$}\},\quad p,q\in S.
\]
Likewise, if $K\subset S$ is a relatively compact subset we define
\[
	\dist_X(p,K):=\inf\{ \dist_X(p,q)\colon q\in K \},\quad p\in S.
\]

An immersed open surface $X\colon S\to\r^n$ $(n\ge 3)$ is said to be {\em complete} if the image by $X$ of any proper path $\gamma\colon [0,1)\to S$ has infinite Euclidean length; this is equivalent to the Riemannian metric induced on $S$ by the Euclidean metric of $\r^n$ via $X$ be complete.

Let $M$ be an open Riemann surface and $n\ge 3$ be an integer. A conformal (i.e. angle-preserving) immersion $X=(X_1,\ldots,X_n)\colon M\to\r^n$ is {\em minimal} (i.e., $X$ has everywhere vanishing mean curvature vector) if, and only if, $X$ is a harmonic map in the classical sense: $\triangle X=0$. Denoting by $\di$ the $\c$-linear part of the exterior differential $d=\di+\overline\di$ on $M$ 
(here $\overline\di$ is the $\c$-antilinear part of $d$), the $1$-form $\di X=(\di X_1,\ldots,\di X_n)$ with values in $\c^n$ is holomorphic, has no zeros, and satisfies
$\sum_{j=1}^n (\di X_j)^2=0$ everywhere on $M$. It follows that the real part $\Re(\di X)$ is an exact real $1$-form on $M$. On the other hand, the {\em flux map} (or simply the {\em flux}) of $X$ is defined as the group homomorphism 
\[
	\Flux_X\colon H_1(M;\z)\to\r^n,
\]
of the first homology group $H_1(M;\z)$ of $M$ with integer coefficients, given by
\[
	\Flux_X(\gamma)=\int_\gamma\Im(\di X)=
	-\imath\int_\gamma \di X,\quad \gamma\in H_1(M;\z),
\]
where $\Im$ denotes the imaginary part and $\imath:=\sqrt{-1}$. We refer e.\ g.\ to Osserman's monograph \cite{Osserman-book} for a standard reference on Minimal Surface Theory.

A compact subset $K\subset M$ is said to be {\em Runge} (also called {\em holomorphically convex} or {\em $\Oscr(M)$-convex}) if its complement $M\setminus K$ has no relatively compact connected components on $M$; by the Runge-Mergelyan theorem \cite{Runge1885AM,Mergelyan1951DAN,Bishop1958PJM} this is equivalent to that every continuous function $K\to\c$, holomorphic in the interior $\mathring K$, may be approximated uniformly on $K$ by holomorphic functions $M\to\c$.

A {\em compact bordered Riemann surface} is a compact Riemann surface $\cM$ with nonempty boundary $bM\subset\overline M$ consisting of finitely many pairwise disjoint smooth Jordan curves; its interior $M=\cM\setminus bM$ is called a {\em bordered Riemann surface}. 
It is classical that every compact bordered Riemann surface $\overline M$ is diffeomorphic to a smoothly bounded compact domain in an open Riemann surface. By a {\em conformal minimal immersion of class $\Cscr^1(\overline M)$} of a compact bordered Riemann surface $\overline M$ into $\r^n$,  we mean a map $\overline M\to\r^n$ of class $\Cscr^1(\overline M)$ whose restriction to $M$ is a conformal minimal immersion. 



\section{Proof of the main results}\label{sec:MR}

In this section we prove Theorems \ref{th:intro-main-v2} and \ref{th:intro-main} in the introduction; both will follow from a recursive application of the following approximation result.

\begin{lemma}\label{lem:main}
Let $D\subset\r^n$ $(n\ge 3)$ be a domain, $\cM=M\cup bM$ be a compact bordered Riemann surface, and  $X\colon \cM\to \r^n$ be a conformal minimal immersion of class $\Cscr^1(\cM)$ such that
\[
	X(\cM)\subset D.
\]
Given a compact domain $K\subset M$, points $p_0\in \mathring K$ and $x_1,\ldots,x_k\in D$ $(k\in\n)$, and numbers $\epsilon>0$ and $\lambda>0$, there is a conformal minimal immersion $Y\colon\cM\to \r^n$ of class $\Cscr^1(\cM)$ satisfying the following conditions:
\begin{enumerate}[\rm (i)]
\item $Y(\cM)\subset D$.
\item $\|Y-X\|_{1,K}<\epsilon$.
\item $\dist(x_j,Y(\cM))<\epsilon\ $ for all $j\in\{1,\ldots,k\}$.
\item $\Flux_Y=\Flux_X$.
\item $\dist_Y(p_0,bM)>\lambda$.
\end{enumerate}
\end{lemma}

We will prove Lemma \ref{lem:main} later in Subsec.\ \ref{ss:Lemma}; we first proceed with the proof of the main results of the paper.


\subsection{Proof of Theorem \ref{th:intro-main-v2} assuming Lemma \ref{lem:main}}

Let $D\subset\r^n$, $M$, $\pgot\colon H_1(M;\z)\to\r^n$, $K\subset M$, $X\colon K\to\r^n$, and $\epsilon>0$ be as in the statement of Theorem \ref{th:intro-main-v2}.

Set $M_0:=K$ and choose an exhaustion of $M$ by connected Runge compact domains $\{M_j \}_{j\in\n}$ such that the Euler characteristic $\chi(M_j\setminus\mathring M_{j-1})\in\{-1,0\}$ for all $j\in\n$ and
\begin{equation}\label{eq:exhaustion}
M_0 \Subset M_1\Subset \cdots \Subset \bigcup_{j\in\z_+}M_j=M.
\end{equation}
Existence of such is well-known; see for instance \cite[Lemma 4.2]{AlarconLopez2013JGA} for a simple proof.

Fix a countable subset $C=\{z_j \}_{j\in\n}\subset D$ with 
\begin{equation}\label{eq:Cdense}
	\overline C=\overline D.
\end{equation}
Set $N_0:=M_0=K$, $Y_0:=X$, and if $n\ge 5$ assume without loss of generality that $Y_0$ is an embedding (as we may in view of \cite[Theorem 1.1]{AlarconForstnericLopez2016MZ}). Also fix a point $p_0\in \mathring N_0$. 

Take a sequence of positive real numbers $\{\epsilon_j\}_{j\in\n}\searrow 0$ which will be specified later. 

We shall recursively construct a sequence $\{N_j,Y_j\}_{j\in\n}$ of smoothly bounded Runge compact domains $N_j\subset M$ and conformal minimal immersions $Y_j\colon N_j\to \r^n$ of class $\Cscr^1(N_j)$ satisfying the following properties for all $j\in\n$:
\begin{enumerate}[\rm (a$_j$)]
	\item $Y_j(N_j)\subset D$.
	\item $N_j\subset M_j$ and $N_j$ is a strong deformation retract of $M_j$.
	\item $\|Y_j-Y_{j-1}\|_{1,N_{j-1}}<\epsilon_j$.
	\item $\dist(z_k,Y_j(N_j))<\epsilon_j$ for all $k\in\{1,\ldots,j\}$.
	\item $\dist_{Y_j}(p_0,bN_j)>j$.
	\item $\Flux_{Y_{j}}(\gamma)=\pgot(\gamma)$ for all closed curves
	 $\gamma\subset N_j$.
	\item If $D=\r^n$ then $N_j= M_j$.
	\item If $n\ge 5$ then $Y_j$ is an embedding.
\end{enumerate}
Observe that condition {\rm (a$_j$)} always holds in case $ D=\r^n$.

Assume for a moment that we have already constructed such a sequence and let us show that if each $\epsilon_j>0$ in the recursive procedure is chosen sufficiently small (in terms of the geometry of $Y_{j-1}$) then the sequence $\{Y_j\}_{j\in\n}$ converges uniformly on compact subsets of
\begin{equation}\label{eq:Omega}
	\Omega:=\bigcup_{j\in\n}N_j \subset M
\end{equation}
to a conformal minimal immersion
\[
	Y:=\lim_{j\to+\infty}Y_j\colon \Omega\to \r^n
\]
satisfying the conclusion of the theorem. Indeed, first of all notice that properties {\rm (b$_j$)}, {\rm (g$_j$)}, \eqref{eq:exhaustion}, and \eqref{eq:Omega} ensure condition {\rm (I)} in the statement of the theorem and that $\Omega=M$ if $D=\r^n$. Now, choosing the $\epsilon_j$'s such that
\begin{equation}\label{eq:epsilon1}
	\sum_{j\in\n}\epsilon_j<\epsilon
\end{equation}
we have in view of {\rm (c$_j$)} that the limit map $Y$ exists and satisfies condition {\rm (II)}. Furthermore, if the sequence $\{\epsilon_j\}_{j\in\n}$ decreases to zero fast enough then, by Harnack's theorem, $Y$ is a conformal minimal immersion. Likewise, by {\rm (c$_j$)}, {\rm (e$_j$)}, and {\rm (f$_j$)}, we have that $Y$ is complete and satisfies {\rm (III)} whenever that each $\epsilon_j>0$ is small enough. 

Let us now check condition {\rm (IV)}. For the first part observe that properties {\rm (a$_j$)} ensure that $Y(\Omega)\subset \overline D$; let us show that $Y(\Omega)\cap bD=\emptyset$. For that, we choose 
\begin{equation}\label{eq:epjYonD}
	\epsilon_j < \frac{1}{j^2} \dist(Y_{j-1}(N_{j-1}),bD)\quad\text{for all $j\in\n$}.  
\end{equation}
Notice that the term in the right hand of the above inequality is positive due to {\rm (a$_j$)}, and hence such an $\epsilon_j>0$ exists. Take $p\in\Omega$ and let us show that $\dist(Y(p),bD)>0$; this will ensure that $Y(\Omega)\subset D$. Choose $j_0\in\n$ such that $p\in N_{j-1}$ for all $j\ge j_0$. Then
\begin{eqnarray*}
	\dist(Y_{j-1}(p),bD) & \le & |Y_{j-1}(p)-Y_j(p) | + \dist(Y_{j}(p),bD)  \\
	&  \stackrel{\textrm{(c$_j$)} }{<} & \epsilon_j + \dist(Y_{j}(p),bD)  \\
	& \stackrel{\eqref{eq:epjYonD}}{<} & \frac{1}{j^2} \dist(Y_{j-1}(p),bD) + \dist(Y_{j}(p),bD).
\end{eqnarray*}
Thus,
$\dist(Y_{j}(p),bD) \ge (1-1/{j^2})\dist(Y_{j-1}(p),bD)$ for all $j\ge j_0$, 
and so
\[
	\dist(Y_{j_0+i}(p),bD) \ge \dist(Y_{j_0}(p),bD) 
	\prod_{j=j_0+1}^{j_0+i}\left(1-\frac{1}{j^2}\right) \quad \text{for all $i\in\n$}.
\]
Taking limits in the above inequality as $i\to+\infty$ we obtain 
\[
	\dist(Y(p),bD)\ge \frac12 \dist(Y_{j_0}(p),bD)>0,
\]
where the latter inequality is ensured by {\rm (a$_{j_0}$)}; take into account that $Y_{j_0}(N_{j_0})$ is compact. This shows that $Y(\Omega)\subset D$.

In order to check the second part of condition {\rm (IV)} pick a point $z\in \overline D$ and a positive number $\delta>0$ and let us show that $\dist(z,Y(\Omega))<\delta$; this will imply that $\overline{Y(\Omega)}=\overline D$. Indeed, in view of \eqref{eq:Cdense} there exists $j_0\in\n$ such that the point $z_{j_0}\in C\subset D$ meets
\begin{equation}\label{eq:zk2}
	|z_{j_0}-z|<\delta/3.
\end{equation}
Moreover, since $\{\epsilon_j\}\searrow 0$, there exists $j_1\in\n$ such that $\epsilon_{j_1}<\delta/3$, and so, for any $j\geq j_1$, {\rm (d$_j$)} guarantees that
\begin{equation}\label{eq:Yzk2}
\dist(z_k,Y_j(N_j))<\delta/3 \quad\text{for all $k\leq j$}.
\end{equation}
Finally, \eqref{eq:epsilon1} ensures the existence of $j_2\in \n$ such that $\sum_{k=j_2}^{\infty}\epsilon_k<\delta/3$ and hence, for all $j > j_2$, properties {\rm (c$_j$)} imply that
\begin{equation}\label{eq:YYj2}
	\| Y-Y_j\|_{1,N_j} <\delta/3.
\end{equation} 
Combining \eqref{eq:zk2}, \eqref{eq:Yzk2}, and \eqref{eq:YYj2} we obtain that, for any $j> \max\{j_0,j_1,j_2\}$,
\begin{eqnarray*}
	\dist(z,Y(\Omega)) & \le & |z-z_{j_0}| + \dist(z_{j_0},Y(\Omega)) \\
	& \stackrel{\eqref{eq:Omega}}{\le} & |z-z_{j_0}| + \dist(z_{j_0},Y(N_j)) \\
	& \le & |z-z_{j_0}| + \dist(z_{j_0},Y_j(N_j))+ \|Y_j-Y\|_{1,N_j}  < \delta.
\end{eqnarray*}
This proves that $Y(\Omega)$ is dense on $\overline D$ and hence condition {\rm (IV)}.

Finally, assume that $n\ge 5$ and let us show that the limit map $Y\colon\Omega\to\r^n$ is one-to-one provided that the positive numbers $\{\epsilon_j\}_{j\in\n}$ are taken sufficiently small. For that it suffices to choose
\begin{equation}\label{eq:epj}
	\epsilon_j < \frac{1}{2j^2} \inf\big\{ |Y_{j-1}(p)-Y_{j-1}(q)| 
	\colon p,q\in N_{j-1},\; {\sf d}(p,q)>1/j \big\}
\end{equation} 
where ${\sf d}(\cdot,\cdot)$ is any fixed Riemannian distance on $M$. Indeed, pick points $p,q\in\Omega$, $p\neq q$, and let us check that $Y(p)\neq Y(q)$. Choose $j_0\in\n$ large enough so that $p,q\in N_{j-1}$ and ${\sf d}(p,q)>1/j$ for all $j\ge j_0$; such exists in view of {\rm (b$_j$)} and \eqref{eq:Omega}. Then
\begin{eqnarray*}
|Y_{j-1}(p)-Y_{j-1}(q) | & \le & |Y_{j-1}(p)-Y_j(p) | + |Y_j(p)-Y_j(q) | + |Y_j(q)-Y_{j-1}(q) | \\
&  \stackrel{\textrm{(c$_j$)} }{<} & 2\epsilon_j + |Y_{j}(p)-Y_{j}(q) |  \\
& \stackrel{\eqref{eq:epj}}{<} & \frac{1}{j^2} |Y_{j-1}(p)-Y_{j-1}(q)| + |Y_{j}(p)-Y_{j}(q) |.
\end{eqnarray*}
As above, this gives that 
$|Y_{j}(p)-Y_{j}(q) | \ge (1-1/{j^2})|Y_{j-1}(p)-Y_{j-1}(q) |$ for all $j\ge j_0$, 
and hence
\[
	|Y_{j_0+i}(p)-Y_{j_0+i}(q) | \ge |Y_{j_0}(p)-Y_{j_0}(q) | \prod_{j=j_0+1}^{j_0+i}\left(1-\frac{1}{j^2}\right) \quad \text{for all $i\in\n$}.
\]
Taking limits we obtain that
\[
	|Y(p)-Y(q)|\ge\frac12 |Y_{j_0}(p)-Y_{j_0}(q)|>0,
\]
where the latter inequality follows from {\rm (h$_{j_0}$)}. This implies that $Y$ is one-to-one, proving condition {\rm (V)} in the statement of the theorem.

To complete the proof it remains to construct the sequence $\{N_j,Y_j\}_{j\in\n}$ satisfying the required properties. We proceed in a recursive way. The basis of the induction is given by the pair $(N_0,Y_0)$ which clearly meets properties {\rm (a$_0$)}, {\rm (b$_0$)}, {\rm (e$_0$)}, {\rm (f$_0$)}, {\rm (g$_0$)}, and {\rm (h$_0$)}; whereas {\rm (c$_0$)} and {\rm (d$_0$)} are vacuous. For the inductive step assume that we have $(N_{j-1},Y_{j-1})$ satisfying {\rm (a$_{j-1}$)}--{\rm (h$_{j-1}$)} and let us construct $(N_j,Y_j)$ enjoying the corresponding properties. We distinguish two different cases depending on the Euler characteristic of $M_j\setminus \mathring  M_{j-1}$.

\smallskip

\noindent \textit{Noncritical case: Assume that $\chi(M_j\setminus\mathring  M_{j-1})=0$}. By the Mergelyan theorem for conformal minimal immersions (see \cite[Theorem 5.3]{AlarconForstnericLopez2016MZ}) we may assume without loss of generality that $Y_{j-1}$ extends, with the same name, to a conformal minimal immersion $M\to\r^n$ with 
\begin{equation}\label{eq:Fluxj-1}
	\Flux_{Y_{j-1}}=\pgot.
\end{equation} 
Next, we choose $N_j\subset M_j$ as any smoothly bounded compact neighborhood of $N_{j-1}$ such that 
\begin{equation}\label{eq:NjD}
	Y_{j-1}(N_j)\subset D
\end{equation}
and that $N_{j-1}$ is a strong deformation retract of $N_j$; such exists in view of {\rm (a$_{j-1}$)}. Since $\chi(M_j\setminus\mathring M_{j-1})=0$, it follows that $N_j$ is a strong deformation retract of $M_j$ as well. This proves {\rm (b$_j$)}. If $D=\r^n$ then we choose, as we may since \eqref{eq:NjD} is always satisfied, $N_j=M_j$, ensuring condition {\rm (g$_j$)}. 

Now, in view of \eqref{eq:NjD}, we may apply Lemma \ref{lem:main} to the domain $D$, the compact bordered Riemann surface $N_j$, the conformal minimal immersion $Y_{j-1}\colon N_j\to D\subset\r^n$ of class $\Cscr^1(N_j)$, the compact domain $N_{j-1}\subset\mathring N_j$, the points $p_0\in \mathring K\subset \mathring N_{j-1}$ and $z_1,\ldots,z_j\in D$, and the positive numbers $\epsilon_j$ and $j>0$.
This provides a conformal minimal immersion $Y_j\colon N_j\to\r^n$ of class $\Cscr^1(N_j)$ enjoying the following properties:
\begin{enumerate}[\rm (i)]
\item $Y_j(N_j)\subset D$.
\item $\|Y_j-Y_{j-1}\|_{1,N_{j-1}}<\epsilon_{j}$.
\item $\dist(z_k,Y_j(N_j))<\epsilon_j$ for all $k\in\{1,\ldots,j\} $.
\item $\Flux_{Y_j}(\gamma)=\Flux_{Y_{j-1}}(\gamma)$ for all closed curves $\gamma\subset N_{j}$.
\item $\dist_{Y_j}(p_0,bN_j)>j$.
\end{enumerate}
Furthermore, we may assume by \cite[Theorem 1.1]{AlarconForstnericLopez2016MZ} that
\begin{enumerate}[\rm (vi)]
\item[\rm (vi)] if $n\ge 5$ then $Y_j$ is an embedding.
\end{enumerate}

We claim that $(N_j,Y_j)$ meets conditions {\rm (a$_j$)}--{\rm (h$_j$)}. Indeed, {\rm (b$_j$)} and {\rm (g$_j$)} are already ensured. On the other hand, conditions {\rm (a$_j$)}, {\rm (c$_j$)}, {\rm (d$_j$)}, {\rm (e$_j$)},  and {\rm (h$_j$)} equal {\rm (i)}, {\rm (ii)}, {\rm (iii)}, {\rm (v)}, and {\rm (vi)}, respectively, whereas {\rm (f$_j$)} is implied by {\rm (iv)} and \eqref{eq:Fluxj-1}. This concludes the proof of the inductive step in the noncritical case.

\smallskip

\noindent \textit{Critical case: Assume that $\chi(M_{j}\setminus\mathring M_{j-1})=-1$}. In this case there is a smooth Jordan arc  $\alpha\subset \mathring M_{j}\setminus \mathring N_{j-1}$, with its two endpoints in $bN_{j-1}$ and otherwise disjoint from $N_{j-1}$, such that
\[
	S:=N_{j-1}\cup\alpha\subset \mathring M_j
\]
is a Runge {\em admissible} subset in $M$ in the sense of \cite[Def.\ 5.1]{AlarconForstnericLopez2016MZ} and a strong deformation retract of $M_j$. 
Fix a nowhere vanishing holomorphic $1$-form $\theta$ on $M$ (such always exists by the Oka-Grauert principle (see \cite[Theorem 5.3.1]{Forstneric2011-book}); for an alternative proof see \cite[Proof of Theorem 4.2]{AlarconFernandezLopez2012CMH}). Next, consider a {\em generalized conformal minimal immersion $(\wt Y,f\theta)$ on $S$} in the sense of \cite[Def.\ 5.2]{AlarconForstnericLopez2016MZ} such that
\[
	\wt Y|_{N_{j-1}}=Y_{j-1},\quad \wt Y(\alpha)\subset D,\quad
	\text{and}\quad  
	\int_\gamma f\theta=\imath\pgot(\gamma) 
	\quad \text{for all closed curves $\gamma$ in $S$.}
\]
Such trivially exists in view of {\rm (a$_{j-1}$)}, {\rm (f$_{j-1}$)}, and the path-connectedness of $D$. By \cite[Theorem 5.3]{AlarconForstnericLopez2016MZ} we may approximate $\wt Y$ in the $\Cscr^1(S)$-topology by conformal minimal immersions $\wt Y_{j-1}\colon M\to\r^n$ having $\pgot$ as flux map and being embeddings if $n\ge 5$. For any close enough such approximation $\wt Y_{j-1}$ of $\wt Y$ there is a compact neighborhood $N_{j-1}'$ of $S$ in $\mathring M_j$ such that $N_{j-1}'\subset M$ is a smoothly bounded Runge compact domain, $S$ is a strong deformation retract of $N_{j-1}'$, and $\wt Y_{j-1}$ formally meets conditions {\rm (a$_{j-1}$)}--{\rm (h$_{j-1}$)} besides {\rm (g$_{j-1}$)}. It follows that the Euler characteristic $\chi(M_j\setminus \mathring N_{j-1}')=0$, which reduces the proof of the inductive step to the noncritical case.

This concludes the recursive construction of the sequence $\{N_j,Y_j\}_{j\in\n}$ with the desired properties, and hence the proof of the theorem.
	
	
\subsection{Proof of Theorem \ref{th:intro-main} assuming Lemma \ref{lem:main}}

Let $K_0\subset M$ be a smoothly bounded compact subset and let $\epsilon>0$. To prove the theorem it suffices to find a complete conformal minimal immersion $Y\colon M\to \r^n$ such that the following conditions are satisfied:
\begin{enumerate}[\rm (a)]
\item $\|Y-X\|_{1,K_0}<\epsilon$.
\item $\Flux_Y=\Flux_X$.
\item $Y(M)\subset D$ and $\overline{Y(M)}=\overline D$. 
\item If $n\ge 5$ then $Y$ is one-to-one.
\end{enumerate}

Up to enlarging $K_0$ if necessary we may assume that $K_0$ is a strong deformation retract of $\overline M$. Pick any countable subset $C=\{z_j\}_{j\in\n}$  of $D$ such that
\begin{equation}\label{eq:Cdense2}
	\overline{C}=\overline{D}.
\end{equation}
Fix a point $p_0\in \mathring K_0\neq \emptyset$ and choose a sequence of positive numbers $\{\epsilon_j\}_{j\in\n}\searrow 0$ that will be specified later. Set $Y_0:=X\colon \cM\to D\subset\r^n$ and, if $n\ge 5$, assume without loss of generality that $Y_0$ is an embedding (cf.\ \cite[Theorem 1.1]{AlarconForstnericLopez2016MZ}). We shall inductively construct a sequence $\{K_j,Y_j\}_{j\in\n}$ of smoothly bounded compact domains 
\begin{equation}\label{eq:exhaustion2}
	K_0\Subset K_1\Subset K_2\Subset\cdots\Subset \bigcup_{j\in\n}K_j=M
\end{equation}
and conformal minimal immersions $\{ Y_j\colon \cM\to \r^n \}_{j\in\n}$ of class $\Cscr^1(\cM)$, satisfying the following properties for all $j\in\n$:
\begin{enumerate}[\rm (I$_j$)]
\item $Y_j(\cM)\subset D$.
\item $\| Y_j-Y_{j-1} \|_{1,K_{j-1}}<\epsilon_j$.
\item $\dist(z_k,Y_j(K_j))<\epsilon_j$ for all $k\in\{1,\ldots, j\}$.
\item $\Flux_{Y_j}(\gamma)=\Flux_{Y_{j-1}}(\gamma)$ for all closed curves $\gamma\subset M$.
\item $\dist_{Y_j}(p_0,bK_j)> j$.
\item If $n\geq 5$ then $Y_j$ is an embedding.
\end{enumerate}

We construct the sequence in an inductive procedure similar to the one in the proof of Theorem \ref{th:intro-main-v2}. The basis of the induction is accomplished by the pair $(K_0,Y_0)$ which clearly satisfies {\rm (I$_0$)}, {\rm (V$_0$)}, and {\rm (VI$_0$)}; conditions {\rm (II$_0$)}, {\rm (III$_0$)}, and {\rm (IV$_0$)} are vacuous. For the inductive step we assume that we already have $(K_{j-1},Y_{j-1})$ satisfying {\rm (I$_{j-1}$)}--{\rm (V$_{j-1}$)}.
By {\rm (I$_{j-1})$ we may apply Lemma \ref{lem:main} to the conformal minimal immersion $Y_{j-1}$, the compact domain $K_{j-1}$, the point $p_0\in\mathring K_0\subset\mathring K_{j-1}$, the points $z_1,\ldots,z_j\in D$, and the positive numbers $\epsilon_j>0$ and $j>0$, obtaining a conformal minimal immersion $Y_j\colon\cM\to\r^n$ of class $\Cscr^1(\cM)$ satisfying the following properties:
\begin{enumerate}[\rm (i)]
	\item $Y_j(\cM)\subset D$.
	\item $\|Y_j-Y_{j-1}\|_{1,K_{j-1}}<\epsilon_j$.
	\item $\dist(z_k,Y(\cM))<\epsilon_j\ $ for all $k\in\{1,\ldots,j\}$.
	\item $\Flux_{Y_j}=\Flux_{Y_{j-1}}$.
	\item $\dist_{Y_j}(p_0,bM)>j$.
\end{enumerate}
Further, by \cite[Theorem 1.1]{AlarconForstnericLopez2016MZ} we may assume that
\begin{enumerate}[\rm (vi)]
\item[\rm (vi)] if $n\ge 5$ then $Y_j$ is an embedding.
\end{enumerate}

Conditions {\rm (I$_{j}$)}, {\rm (II$_j$)}, {\rm (IV$_{j}$)}, and {\rm (VI$_{j}$)} equal {\rm (i)}, {\rm (ii)}, {\rm (iv)}, and {\rm (vi)}. Finally, since the inequalities in {\rm (iii)} and {\rm (v)} are both strict, conditions {\rm (III$_j$)} and {\rm (V$_j$)} hold for any large enough smoothly bounded compact domain $K_j\subset M$ being a strong deformation retract of $\cM$. At each step in the recursive construction, we choose such a $K_j$ containing $K_{j-1}$ in its interior and being large enough so that \eqref{eq:exhaustion2} is satisfied. This closes the induction and  concludes the construction of the sequence $\{K_j,Y_j\}_{j\in\n}$ satisfying conditions {\rm (I$_j$)}--{\rm (VI$_j$)}.

We claim that choosing the number $\epsilon_j>0$ sufficiently small (depending on the geometry of $Y_{j-1}$) at each step in the recursive construction, the sequence $\{Y_j\}_{j\in\n}$ converges uniformly on compact subsets in $M$ to a limit map
\[
     Y:=\lim\limits_{j\to\infty} Y_j\colon M\to\r^n
\]
which satisfies conditions {\rm (a)}--{\rm(d)}. Indeed, reasoning as in the proof of Theorem \ref{th:intro-main-v2}, {\rm (II$_j$)} ensures that the limit map $Y$ is a conformal minimal immersion and meets {\rm (a)}. 
On the other hand, {\rm (IV$_j$)} implies {\rm (b)}; {\rm (V$_j$)} and {\rm (II$_j$)} guarantee the completeness of $Y$; {\rm (c)} follows from  {\rm (I$_j$)},  {\rm (II$_j$)}, and {\rm (III$_j$)}; and properties {\rm (II$_j$)} and {\rm (VI$_j$)} give condition {\rm (d)}. This completes the proof.


\subsection{Proof of Lemma \ref{lem:main}}\label{ss:Lemma}

Without loss of generality we may assume that $k=1$; the general case follows from a standard finite recursive application of this particular one. Call $x:=x_1$

We may also assume without loss of generality that $\cM$ is a smoothly bounded compact domain in an open Riemann surface $\Rcal$. Pick a point $p\in bM$ and  a smooth embedded arc $\gamma\subset \Rcal\setminus M$ having $p$ as an endpoint, being otherwise disjoint from $\cM$, and such that 
\[
	S:=\overline M\cup\gamma
\]
is a Runge {\em admissible} subset of $\Rcal$ in the sense of \cite[Def.\ 5.1]{AlarconForstnericLopez2016MZ}. Let $q\in\Rcal\setminus\cM$ denote the other endpoint of $\gamma$.
	
Fix a nowhere vanishing holomorphic $1$-form $\theta$ on $\Rcal$. Consider a {\em generalized conformal minimal immersion $(\wt X,f\theta)$ on $S$} in the sense of \cite[Def.\ 5.2]{AlarconForstnericLopez2016MZ} such that the $\Cscr^1(S)$-map $\wt X \colon S\to\r^n$ satisfies the following properties: 
\begin{enumerate}[AAA]
\item[\rm (A)] $\wt X|_{\cM}=X$.
\item[\rm (B)] $\wt X|_{\gamma}\subset D$.
\item[\rm (C)] $\wt X(q)=x$.
\end{enumerate}
Existence of such is trivial; recall that $X(\cM)\subset D$ and that $D$ is path-connected.

Fix a constant $\delta>0$ to be specified later.

The Runge-Mergelyan theorem for conformal minimal immersions \cite[Theorem 5.3]{AlarconForstnericLopez2016MZ} provides a conformal minimal immersion $\wt Y\colon \Rcal\to\r^n$ such that
\begin{enumerate}[AAA]
\item[\rm (D)] $\|\wt Y - \wt X\|_{1,S}<\delta$ and
\vspace{1mm}
\item[\rm (E)] $\Flux_{\wt Y}(\alpha)=\Flux_{\wt X}(\alpha)$ for all closed curves $\alpha\subset M$.
\end{enumerate}
Since $X$ assumes values in $D$, properties {\rm (A)} and {\rm (B)} ensure that $\wt X(S)\subset D$, and hence, choosing $\delta>0$ sufficiently small, {\rm (D)} guarantees the existence of a small open neighborhood $U$ of $S$
in $\Rcal$ such that
\begin{equation}\label{eq:YUD}
\wt Y(U)\subset D.
\end{equation} 

Next we use the method of exposing boundary points on a compact bordered Riemann surface. Choose small open neighborhoods $W'\Subset W\Subset U\setminus K$ and $V\Subset U$ of $p$ and $\gamma$  in $U$, respectively. By Forstneri\v c and Wold \cite[Theorem 2.3]{ForstnericWold2009JMPA} (see also \cite[Theorem 8.8.1]{Forstneric2011-book}) there exists a smooth diffeomorphism  
\begin{equation}\label{eq:phi(M)}
	\phi\colon \cM\to\phi(\cM)\subset U
\end{equation}
satisfying the following properties (see Figure \ref{fig:pic1}):
\begin{enumerate}[AAA]
\item[\rm (F)] $\phi\colon M \to\phi(M)$ is a biholomorphism.
\item[\rm (G)] $\phi$ is $\delta$-close to the identity in the $\Cscr^1$-norm on $\cM\setminus W'$.
\item[\rm (H)] $\phi(p)=q\in b\phi(\cM)$ and $\phi(\cM\cap W')\subset W\cup V$.
\end{enumerate}

\begin{figure}[H]
	\includegraphics[height=4.5cm]{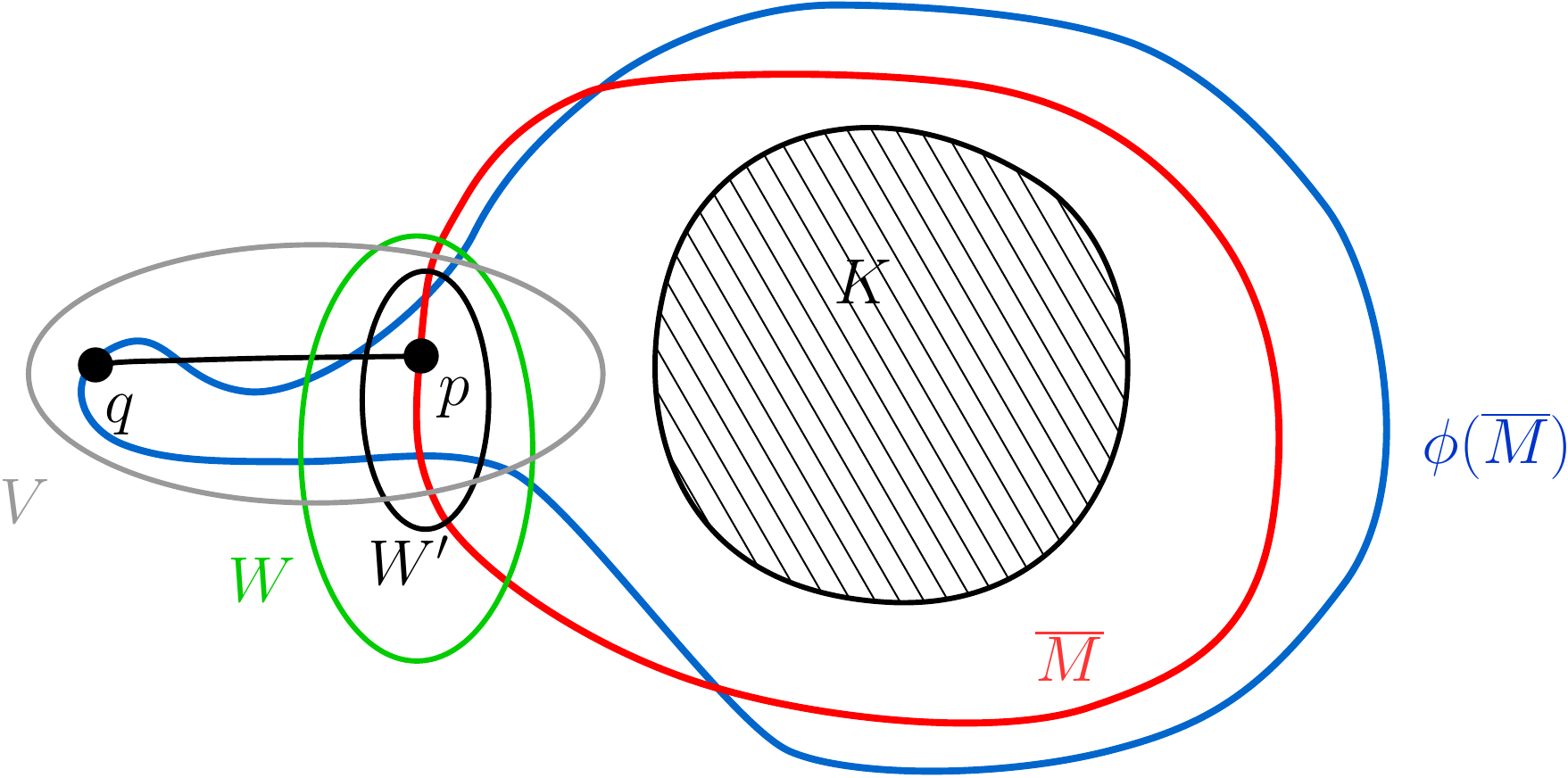}
	\caption{The diffeomorphism $\phi\colon \cM\to\phi(\cM)\subset U$}\label{fig:pic1}
\end{figure}

We claim that the conformal minimal immersion $\wt Y\circ\phi\colon\cM\to\r^n$ of class $\Cscr^1(\cM)$ formally satisfies conditions {\rm (i)}--{\rm (iv)} in the statement of the lemma provided that $\delta>0$ is chosen sufficiently small. Indeed, by \eqref{eq:YUD} and \eqref{eq:phi(M)} we have that $\wt Y(\phi(\cM))\subset D$, proving {\rm (i)}. On the other hand, since $K\subset \cM\setminus W'$, properties {\rm (G)}, {\rm (D)}, and {\rm (A)} give that $\|\wt Y\circ\phi-X\|_{1,K}<\epsilon$, whenever that $\delta>0$ is small enough, which ensures condition {\rm (ii)}. Finally, properties {\rm (H)}, {\rm (D)}, and {\rm (C)} guarantee {\rm (iii)} for any $\delta<\epsilon$, whereas {\rm (F)}, {\rm (E)}, and {\rm (A)} imply {\rm (iv)}.
	
Finally, \cite[Lemma 4.1]{AlarconDrinovecForstnericLopez2015PLMS} enables us to approximate the immersion $\wt Y\circ\phi\colon\cM\to\r^n$ in the $\Cscr^0(\cM)$-topology, and hence in the $\Cscr^1(K)$-topology, by conformal minimal immersions $Y\colon\cM\to\r^n$ of class $\Cscr^1(\cM)$ satisfying {\rm (v)} and $\Flux_Y=\Flux_{\wt Y\circ \phi}$; the latter ensures {\rm (iv)}. It is clear that any close enough such approximation $Y$ of $\wt Y\circ \phi$ still satisfies conditions {\rm (i)}, {\rm (ii)}, and {\rm (iii)}. This concludes the proof of Lemma \ref{lem:main}.

The proof of Theorems \ref{th:intro-main-v2} and \ref{th:intro-main} is now complete.


\section{Analogous results for other families of surfaces}\label{sec:results}

As we already pointed out in the introduction of this paper, all the tools required in the proof of Theorems \ref{th:intro-main-v2} and  \ref{th:intro-main} (i.e. the Runge-Mergelyan approximation, the general position result, and the Riemann-Hilbert method) are also available for some other interesting objects; namely, non-orientable minimal surfaces, complex curves, and holomorphic null and Legendrian curves. Therefore, our techniques easily adapt to give analogous results to Theorems \ref{th:intro-main-v2} and  \ref{th:intro-main} for all these families of surfaces; we shall now discuss some of them, leaving the details of the proofs to the interested reader.
%
%

\subsection{Non-orientable minimal surfaces in $\r^n$} These surfaces appeared in the very origin of Minimal Surface Theory (we refer to the seminal paper by Lie \cite{Lie1878MA} from 1878) and there is a large literature devoted to their study. {\em Conformal non-orientable minimal surfaces} in $\r^n$ for $n\ge 3$ are characterized as the images of conformal minimal inmmersions $X\colon M\to\r^n$ such that $X\circ\Igot=X$, where $\Igot\colon M\to M$ is an antiholomorphic involution without fixed points on an open Riemann surface $M$. For such an immersion we have that
\begin{equation}\label{eq:Flux-no}
	\Flux_X(\Igot_*\gamma)=-\Flux_X(\gamma)\quad \text{for all $\gamma\in H_1(M;\z)$}.
\end{equation}

Rencently, Alarc\'on, Forstneri\v c, and L\'opez introduced in \cite{AlarconForstnericLopez2016Pre1} new complex analytic techniques in the study of non-orientable minimal surfaces in $\r^n$; in particular, they provided all the required tools in our method of proof (see also \cite{AlarconLopez2015GT} for the Runge-Mergelyan approximation in dimension $3$). As happens in the orientable case, the general position of non-orientable minimal surfaces is embedded in $\r^n$ for all $n\ge 5$. Thus, completely analogous results to Theorems \ref{th:intro-main-v2} and  \ref{th:intro-main} may be proved in the non-orientable framework under the necessary condition \eqref{eq:Flux-no} on the flux map.

%
%
\subsection{Complex curves in $\c^n$} All the above mentioned tools are classical for holomorphic immersions of open Riemann surfaces into $\c^n$ for $n\ge 2$, being embedded the general position for $n\ge 3$. We refer to Bishop \cite{Bishop1958PJM} for the Runge-Mergelyan approximation (see also \cite{Runge1885AM,Mergelyan1951DAN}) and to \cite{DrinovecForstneric2007DMJ,AlarconForstneric2013MA,AlarconForstneric2015AS} for the Riemann-Hilbert method (see also the introduction of Drinovec Drnov\v sek and Forstneri\v c \cite{DrinovecForstneric2012IUMJ} for a survey on this subject).

For example, by following the proof of Theorem \ref{th:intro-main-v2} one may show the following
\begin{theorem}\label{th:intro-complex}
	Let $M$ be an open Riemann surface. The set of complete holomorphic immersions $M\to\c^n$ $(n\ge 2)$ with dense images forms a dense subset in the set $\Oscr(M,\c^n)$ of all holomorphic functions $M\to\c^n$ endowed with the compact-open topology. Furthermore, if $n\ge 3$ then the set of all complete holomorphic one-to-one immersions $M\to\c^n$ with dense images is also dense in $\Oscr(M,\c^n)$. 
\end{theorem}
We emphasize that the novelty of Theorem \ref{th:intro-complex} is that it concerns {\em complete} immersions; obviously, the set of all holomorphic immersions $M\to\c^n$ is much larger than the subset consisting of the complete ones. Indeed, without completeness, there are many general such results in the literature. For instance, if we consider the space $\Oscr(S,Z)$ of all holomorphic maps of a {\em Stein manifold} $S$ (we refer to Gunning and Rossi \cite{GunningRossi2009AMS} and H\"ormander \cite{Hormander1990Book} for the theory of Stein manifolds) into an {\em Oka manifold} $Z$, endowed with the compact-open topology, then the basic Oka property with approximation and interpolation (see \cite[Theorem 5.4.4]{Forstneric2011-book}) easily implies that those maps in $\Oscr(S,Z)$ having dense image form a dense subset; further, if $\dim Z\ge 2\dim S$ (respectively, $\dim Z\ge 2\dim S+1$) then, by general position (see \cite[Theorem 7.9.1 and Corollary 7.9.3]{Forstneric2011-book}), the subset of immersions (respectively, one-to-one immersions) with dense image is also dense in $\Oscr(S,Z)$. On the other hand, if $\dim S\ge \dim Z$ then there are strongly dominating {\em surjective} holomorphic maps $S\to Z$ (see Forstneri\v c \cite{Forstneric2016surjective} and the references therein).

In the same line, Forstneri\v c and Winkelmann proved in \cite{ForstnericWinkelmann2005MRL} that, for any connected complex manifold $Z$ (not necessarily Oka), the set of all holomorphic maps of the unit disk $\d\subset\c$ into $Z$ with dense images is dense in $\Oscr(\d,Z)$; see also Winkelmann \cite{Winkelmann2005MZ} for a previous partial result in this direction. 

%
%

\subsection{Holomorphic null curves in $\c^n$} These are holomorphic immersions $F=(F_1,\ldots,F_n)\colon M\to\c^n$ $(n\ge 3)$ of an open Riemann surface $M$ into $\c^n$ which are directed by the null quadric 
\[
	\Agot=\{z=(z_1,\ldots,z_n)\in\c^n\colon z_1^2+\cdots+z_n^2=0\};
\]
equivalently, satisfying the nullity condition
\[
	(dF_1)^2+\cdots+(dF_n)^2=0\quad \text{everywhere on $M$}.
\]
Notice that the punctured null quadric $\Agot_*=\Agot\setminus\{0\}$ is an Oka manifold (see \cite[Example 4.4]{AlarconForstneric2014IM}). These curves are closely related to minimal surfaces in $\r^n$ since the real and the imaginary part of a null curve $M\to \c^n$ are flux-vanishing conformal minimal immersions $M\to\r^n$ (see e.\ g.\ Osserman \cite{Osserman-book}).
The required tools in order to prove analogous results to Theorems \ref{th:intro-main-v2} and  \ref{th:intro-main} for holomorphic null curves have been provided recently in \cite{AlarconLopez2012JDG,AlarconForstneric2014IM,AlarconForstneric2015MA,AlarconDrinovecForstnericLopez2015PLMS}. In this framework, the general position is embedded for $n\ge 3$.

%
%

\subsection{Holomorphic Legendrian curves in $\c^{2n+1}$} 

These are holomorphic immersions $F=(X_1,Y_1,\ldots,X_n,Y_n,Z)\colon M\to\c^{2n+1}$ $(n\in\n$) of an open Riemann surface $M$ into $\c^{2n+1}$ which are tangent to the standard holomorphic contact structure of $\c^{2n+1}$; equivalently, such that
\[
	dZ+\sum_{j=1}^n X_j\, dY_j=0\quad \text{everywhere on $M$}.
\]
All the needed tools in this case were furnished by Alarc\'on, Forstneri\v c, and L\'opez in
\cite{AlarconForstnericLopez2016Legendrian}, being the general position embedded for all $n\in\n$. Holomorphic Legendrian curves are complex analogues of real Legendrian curves in $\r^{2n+1}$ which play an important role in differential geometry; in particular, in contact geometry.

Recall that a {\em complex contact manifold} is a complex manifold $W$ of odd dimension $2n+1\ge 3$ endowed with a {\em holomorphic contact structure} $\Lscr$; the latter is a holomorphic vector subbundle of complex codimension one in the tangent bundle $TW$ such that every point $p\in W$ admits an open neighborhood $U\subset W$ in which $\Lscr|_U=\ker\eta$ for a holomorphic $1$-form $\eta$ on $U$ satisfying $\eta\wedge (d\eta)^n\ne 0$ everywhere on $U$. A holomorphic immersion $F\colon M\to W$ is said to be {\em Legendrian} if it is everywhere tangent to the contact structure:
\[
	dF_p(T_pM)\subset \Lscr_{F(p)}\quad \text{for all $p\in M$}.
\]
By Darboux's theorem (see \cite[Theorem A.2]{AlarconForstnericLopez2016Legendrian}) every complex contact manifold $(W,\Lscr)$ of dimension $2n+1$ is locally contactomorphic to $\c^{2n+1}$ endowed with its standard holomorphic contact structure. Thus, as a direct consequence of the analogous to Theorems \ref{th:intro-main-v2} and \ref{th:intro-main} for Legendrian curves in $\c^{2n+1}$ one easily obtains the following
\begin{corollary}
Let $(W,\Lscr)$ be a complex contact manifold. Every point $p\in W$ admits an open neighborhood $U\subset W$ with the following property: Given a domain $V\Subset U$ there are holomorphic Legendrian one-to-one immersions $M\to V$ which are dense on $V$ and are complete with respect to every Riemannian metric in $W$, where $M$ is either a given bordered Riemann surface or some complex structure on any given smooth orientable connected open surface.
\end{corollary}
The proof of the above corollary follows the one of \cite[Corollary 1.3]{AlarconForstnericLopez2016Legendrian}; we refer there for the details. It remains as an open question whether every complex contact manifold, endowed with a Riemannian metric, admits complete dense complex Legendrian curves.


\subsection*{Acknowledgements}
A.\ Alarc\'on is supported by the Ram\'on y Cajal program of the Spanish Ministry of Economy and Competitiveness.
A.\ Alarc\'on and I.\ Castro-Infantes are partially supported by the MINECO/FEDER grant no. MTM2014-52368-P, Spain. 

We thank Franc Forstneri\v c, Francisco J. L\'opez, and Joaqu\'in P\'erez for helpful suggestions which led to improvement of the paper.




\medskip

\noindent Antonio Alarc\'{o}n

\noindent Departamento de Geometr\'{\i}a y Topolog\'{\i}a e Instituto de Matem\'aticas (IEMath-GR), Universidad de Granada, Campus de Fuentenueva s/n, E--18071 Granada, Spain.

\noindent  e-mail: {\tt alarcon@ugr.es}

\medskip

\noindent Ildefonso Castro-Infantes

\noindent Departamento de Geometr\'{\i}a y Topolog\'{\i}a, 
Universidad de Granada, Campus de Fuentenueva s/n, E--18071 Granada, Spain.

\noindent  e-mail: {\tt icastroinfantes@ugr.es}

\end{document}